\documentclass[mathpazo]{cicp}

\usepackage{graphicx}
\usepackage{color}
\usepackage{enumerate}
\usepackage{multirow}
\usepackage{tabls}
\usepackage{booktabs}

\begin{document}
\title{A High-Order Modified  Finite Volume WENO
Method on 3D Cartesian Grids}

 \author[Du Y L et.~al.]{Yulong Du\affil{1},
       Li Yuan\affil{2}\comma\affil{3}\comma\corrauth, Yahui Wang\affil{2}\comma\affil{3}}
 \address{\affilnum{1}\ School of Mathematics and Systems Science,
           Beihang University,
           Beijing  100191, P. R. China.\\
          \affilnum{2}\ ICMSEC and LSEC,
          Academy of Mathematics and Systems Science,
          Chinese Academy of Sciences,
          Beijing 100190, P. R. China.\\
          \affilnum{3}\ School of Mathematical Sciences,
          University of Chinese Academy of Sciences,
          Beijing 100190, P. R. China.}


 \emails{{\tt kunyu0918@163.com} (Yulong~Du),
          {\tt lyuan@lsec.cc.ac.cn} (Li~Yuan),  {\tt wangyh14@lsec.cc.ac.cn} (Yahui~Wang)}

\begin{abstract}
The modified dimension-by-dimension finite volume (FV) WENO method on Cartesian grids proposed by Buchm\"{u}ller and Helzel can retain the full order of accuracy of the one-dimensional WENO reconstruction and requires only one flux computation per interface. The high-order accurate conversion between face-averaged
values and face-center point values is the main ingredient of this method. In this paper, we derive sixth-order accurate conversion formulas on three-dimensional Cartesian grids.  It is shown that the resulting modified FV WENO method is efficient and high-order accurate when applied to smooth nonlinear multidimensional problems, and  is robust for calculating non-smooth nonlinear problems with strong shocks.
\end{abstract}


\maketitle

\section{Introduction}
\label{sec1}
The standard weighted essentially non-oscillatory (WENO) method proposed by Jiang and Shu \cite{shu96} is widely used for solving hyperbolic
conservation laws. The simplest way to use WENO methods on multidimensional Cartesian grids is to apply a one-dimensional WENO scheme in each direction \cite{shu97}. Conservative finite difference WENO methods based on flux interpolation are used in a dimension-by-dimension fashion and they can retain the full order accuracy of the one-dimensional WENO scheme for linear as well as nonlinear multi-dimensional conservation laws. However, in some situations such as adaptively refined Cartesian grids and multi-block Cartesian grids, finite volume methods (FVMs) are more convenient than finite difference methods as FVMs  admit a simple formulation around hanging nodes.  Unfortunately, while FV WENO methods based on a dimension-by-dimension fashion retain the full order of accuracy for smooth solutions of linear multi-dimensional problems,  they are only
second-order accurate for smooth solutions of nonlinear multi-dimensional problems \cite{shu11,improved14}.

A high-order FVM generally  includes variable reconstruction within a cell (k-exact reconstruction \cite{barth90} and its variants \cite{hushu99, ccp1, ccp2, ccp3, ccp4, ccp5}) and high order flux quadrature on the cell interfaces. On Cartesian grids, the expensive  multi-dimensional WENO reconstruction
is not necessary. Instead, a series of one-dimensional WENO reconstructions are applied in all directions in order to obtain high-order accurate point values of the conserved quantities at the quadrature points of a cell interface, and then evaluate numerical fluxes at these quadrature points.  However, the computational cost of such high-order FV WENO methods on Cartesian grids is still large \cite{shu11,yuan2011,toro04}.

Recently, Buchm\"{u}ller and Helzel  \cite{improved14} proposed a modification to the dimension-by-dimension
FV  WENO method on Cartesian grids  and applied this modified method on adaptive Cartesian meshes \cite{improved16, improved18}. Later on a fourth-order quadrature modification flux (QMF) method was introduced and applied on adaptive Cartesian meshes by Tamaki and Imamura \cite{Tamaki2017CAF}. A key technique used in  Refs.~\cite{improved14,improved16, improved18} is the conversion between face-averaged values and face-centered values, which helps improve the spatial order of accuracy of the dimension-by-dimension FV WENO method. However, Refs. \cite{improved14, improved16} mainly concentrated on two-dimensional problems and  Ref. \cite{improved18} only gave the fourth-order conversion formulas  on three-dimensional (3D) Cartesian  grids.
In this paper,  we  further develop the modified FV WENO method by deriving sixth-order conversion formulas on 3D Cartesian grids which are not available in Refs.~\cite{improved14, improved16, improved18}. The derivation is
based on the observation that the differentiation and cell-averaging are exchangeable \cite{Tamaki2017CAF}, and it  can be extended to  even higher order accuracy of conversion.
Furthermore, we use the characteristic variables as the reconstructed quantities for the system of conservation laws. For the temporal discretization we use the same  Runge-Kutta methods of order fifth or seven as Ref.~\cite{improved14}.

The rest of this paper is organized as follows. In Section 2, the modified FV WENO method is explained, and the sixth-order formulas for conversion between face-averaged values and face center point values on 3D Cartesian grids are derived.  Numerical results are presented in Section 3 to verify the accuracy, efficiency and robustness of the modified method. Concluding remarks are given in Section 4.

\section{Modified finite volume WENO method}
\label{sec:2}
In this section we first give the standard dimension-by-dimension
FV WENO method, and then derive sixth-order conversion formulas on 3D Cartesian grids. Finally, we build up our modified dimension-by-dimension FV WENO method.
\subsection{Dimension-by-dimension finite volume WENO method}
\label{sec:2.1}
The 3D system  of conservation laws with initial conditions  are given by
\begin{equation}
\label{se1:eq1}
\begin{array}{l}
\displaystyle\partial_{t}{u}+\partial_{x}f(u)+\partial_{y}g(u)+\partial_{z}h(u)=0,\\
\displaystyle u(x,y,z,0)=u_{0}(x,y,z),
\end{array}
\end{equation}
where $u(x,y,z,t): \mathbb{R}^3\times\mathbb{R}^{+}\rightarrow\mathbb{R}^{m}$ is a vector of conserved quantities, and $f(u)$, $g(u)$, $h(u)$: $\mathbb{R}^{m}\rightarrow\mathbb{R}^{m}$ are vector valued flux functions.

Let $C_{i,j,k}=(x_{i-\frac{1}{2}}, x_{i+\frac{1}{2}})\times (y_{j-\frac{1}{2}}, y_{j+\frac{1}{2}})\times (z_{k-\frac{1}{2}}, z_{k+\frac{1}{2}})$  be a control volume in the $xyz$ space,
 with uniform grid sizes $\Delta x=x_{i+\frac{1}{2}}-x_{i-\frac{1}{2}}$, $\Delta y=y_{j+\frac{1}{2}}-y_{j-\frac{1}{2}}$ and $\Delta z=z_{k+\frac{1}{2}}-z_{k-\frac{1}{2}}$.
Integrating equation (\ref{se1:eq1}) over $C_{i,j,k}$, we obtain a finite volume method in the semi-discrete form
\begin{equation}
\label{se2:eq3}
\begin{split}
\displaystyle\frac{d}{dt}U_{i,j,k}(t)&=\displaystyle-\frac{1}{\Delta x}\left(\hat{F}_{i+\frac{1}{2},j,k}(t)-\hat{F}_{i-\frac{1}{2},j,k}(t) \right) -\frac{1}{\Delta y}\left(\hat{G}_{i,j+\frac{1}{2},k}(t)-\hat{G}_{i,j-\frac{1}{2},k}(t) \right)\\
&\quad -\frac{1}{\Delta z}\left(\hat{H}_{i,j,k+\frac{1}{2}}(t)-\hat{H}_{i,j,k-\frac{1}{2}}(t) \right),
\end{split}
\end{equation}
where $U_{i,j,k}(t)$ is the cell average of the conserved quantities, and $ \hat{F}_{i+1/2,j,k}(t)$, $\hat{G}_{i,j+1/2,k}(t)$ and $\hat{H}_{i,j,k+1/2}(t)$ are face-averaged numerical fluxes. 

In this work, explicit high-order Runge-Kutta methods are used for the temporal discretization of Eq.~(\ref{se2:eq3}). For the spatial discretization we use a one-dimensional piecewise polynomial WENOZ reconstruction \cite{wenoz,wenoz7} in each direction. For example, we construct polynomials $q^{1}_{i,j,k}(x)$ in the $x$ direction, which are local approximations of the $yz$ plane-averaged value of the conserved quantities  $u(x,y,z,t)$ in cell $C_{i,j,k}$.  
By evaluating the polynomials at the interfaces, 
we get two reconstructed face-averaged values of the conserved quantities,
\begin{eqnarray}
\label{se2:eq5}
U^{+}_{i-\frac{1}{2},j,k}:=q^{1}_{i,j,k}(x_{i-\frac{1}{2}}), \quad U^{-}_{i+\frac{1}{2},j,k}:=q^{1}_{i,j,k}(x_{i+\frac{1}{2}}).
\end{eqnarray}

Let $\tilde{u}_{i+1/2,j,k}$ denote the exact face-averaged values of the conserved quantities $u$ at a cell interface
\begin{equation}
\label{exactface}
 \tilde{u}_{i+\frac{1}{2},j,k}=\frac{1}{\Delta y\Delta z}\int_{z_{k-\frac{1}{2}}}^{z_{k+\frac{1}{2}}}
\int_{y_{j-\frac{1}{2}}}^{y_{j+\frac{1}{2}}} u(x_{i+\frac{1}{2}},y,z)dy dz.
\end{equation}
Then the reconstructed face-averaged values satisfy (for sufficiently smooth functions $u$)
\begin{eqnarray}
\label{se2:eq6}
U^{\pm}_{i+\frac{1}{2},j,k}=\tilde{u}_{i+\frac{1}{2},j,k}+\mathcal{O}(\Delta x^{p})
\end{eqnarray}
with $p=7$ in this paper.

The numerical flux can be obtained by using a numerical flux function defined as $\mathcal{F}(u^{-},u^{+})$, which is at least Lipschitz continuous and consistent with the physical flux $f$ in the sense that $\mathcal{F}(u,u)=f(u)$. In this paper we use the Lax-Friedrichs flux or the HLLC flux in respective numerical examples.

It is well known that the dimension-by-dimension FV WENO approach can not retain the $p^{\text{th}}\text{-order}$ accuracy of the underlying reconstruction polynomial. It is only second-order accurate for nonlinear multidimensional conservation laws. Only in the linear case, i.e. $f(u)=Au$ with a constant matrix $A\in\mathbb{R}^{m\times m}$, it can retain the full spatial order of accuracy.
This has been shown for the Lax-Friedrichs flux in Ref.~\cite{improved14} and for a general consistent  numerical flux function in Ref.~\cite{shu11}.

\subsection{Conversion between average values and point values}
\label{sec:2.2}
 Noting that the fourth-order conversion formulas between  cell average values and cell center point values for any function $q(x,y)$ of two-dimensional variables have already been given in Ref.\cite{improved18}, we concentrate on sixth-order conversion formulas. The derived formulas will be used at cell interfaces of a 3D Cartesian grid.

Denote with $Q(x,y)$ an integrated function of the function $q(x,y)$  as
\begin{equation}
\label{djjiaohuan}
\displaystyle Q(x,y):=\overline{q}(x,y)=\frac{1}{\Delta x\Delta y}\int_{y-\frac{\Delta y}{2}}^{y+\frac{\Delta y}{2}}\int_{x-\frac{\Delta x}{2}}^{x+\frac{\Delta x}{2}}  q(\xi,\eta)d\xi  d\eta.
\end{equation}
In a grid cell $(i,j)$, i.e., the rectangle $(x_{i-\frac{1}{2}}, x_{i+\frac{1}{2}})\times(y_{j-\frac{1}{2}}, y_{j+\frac{1}{2}})$, let $Q(x_{i},y_{j})$ denote a cell average value and $q_{i,j}$ be the cell center point value of the function $q(x,y)$. For sufficiently smooth functions $q(x,y): \mathbb{R}^{2}\rightarrow\mathbb{R}^{m}$, Taylor series expansion provides
\begin{equation}
\label{se3:eq4}
\setlength{\arraycolsep}{0.5pt}
\begin{array}{ll}
Q_{i,j}:&=Q(x_{i},y_{j})\displaystyle=\frac{1}{\Delta x\Delta y}\int_{y_{j-\frac{1}{2}}}^{y_{j+\frac{1}{2}}}
\int_{x_{i-\frac{1}{2}}}^{x_{i+\frac{1}{2}}} q(x,y)dxdy\\
&\displaystyle=\frac{1}{\Delta x\Delta y}\int_{-\frac{\Delta y}{2}}^{\frac{\Delta y}{2}}\int_{-\frac{\Delta x}{2}}^{\frac{\Delta x}{2}}q(x_{i}+x,y_{j}+y)dxdy\\
&\displaystyle=\frac{1}{\Delta x\Delta y}{ \int_{-\frac{\Delta y}{2}}^{\frac{\Delta y}{2}}\int_{-\frac{\Delta x}{2}}^{\frac{\Delta x}{2}} }  \left\{q(x_{i},y_{j})+xq_x(x_{i},y_{j})+ yq_y (x_{i},y_{j})\right.\\ &~~~\displaystyle+\frac{1}{2}\left[x^{2}q_{xx} (x_{i},y_{j})+2xy q_{xy}(x_{i},y_{j})+ y^2 q_{yy}  (x_{i},y_{j}) \right]\\
&~~~\displaystyle+\frac{1}{6}\left[x^{3}q_{xxx} (x_{i},y_{j})+ 3x^2y q_{xxy} (x_{i},y_{j})+3xy^2q_{xyy} (x_{i},y_{j})\right.\\
&~~~\displaystyle\left.+y^{3}q_{yyy} (x_{i},y_{j})\right]+\frac{1}{24}\left[x^{4}q_{xxxx}(x_{i},y_{j})+4x^{3}yq_{xxxy}(x_{i},y_{j})\right.\\
&~~~\displaystyle\left.+ 6x^{2}y^{2}q_{xxyy}(x_{i},y_{j})+4xy^{3}q_{xyyy}(x_{i},y_{j})+y^{4}q_{yyyy}(x_{i},y_{j})\right]\\
&~~~\displaystyle\left.+\cdots\right\}dxdy
\end{array}
\end{equation}
and thus the transformation
\begin{equation}
\label{se3:eq5}
\setlength{\arraycolsep}{0.5pt}
\begin{array}{ll}
Q_{i,j}&=\displaystyle q_{i,j}+\frac{\Delta x^2}{24}q_{xx}(x_{i},y_{j})+\frac{\Delta y^2}{24}q_{yy}(x_{i},y_{j})+\frac{\Delta x^4}{1920}q_{xxxx} (x_{i},y_{j})\\
&~~~\displaystyle+\frac{\Delta x^2\Delta y^2}{576} q_{xxyy} (x_{i},y_{j})+\frac{\Delta y^4}{1920}q_{yyyy} (x_{i},y_{j})\\
&~~~\displaystyle+\mathcal{O}(\Delta x^6+\Delta x^4\Delta y^2+\Delta x^2\Delta y^4+\Delta y^6)
\end{array}
\end{equation}
between point values $q_{i,j}$ and cell average values $Q_{i,j}$.

In order to derive sixth-order accurate conversion formulas, we need the approximations of the second and fourth partial derivatives of $q(x,y)$ at the cell center point in Eq.~(\ref{se3:eq5}). It is easy to approximate these derivatives from point values by using standard finite difference schemes directly. However, if we transform from cell average values to point values, it is not trivial to express $q_{xxxx}(x_i,y_j)$ etc. in terms of cell average  values.

Thanks to the notation \cite{Tamaki2017CAF} that the differentiation and cell-averaging  are exchangeable in (\ref{djjiaohuan}), e.g., $Q_{xx}=\overline{q_{xx}},~Q_{xxyy}=\overline{q_{xxyy}}$, we can proceed like  (\ref{se3:eq4}) to derive similar transformations between point values and cell average values of the partial derivatives,
\begin{equation}
\label{d2}
\setlength{\arraycolsep}{0.5pt}
\begin{array}{ll}
\displaystyle Q_{xxi,j}=&\displaystyle q_{xxi,j}+\frac{\Delta x^2}{24}q_{xxxx}(x_{i},y_{j})+\frac{\Delta y^2}{24}q_{xxyy}(x_{i},y_{j})\\
&+\mathcal{O}(\Delta x^4+\Delta x^2\Delta y^2+\Delta y^4),\\
\displaystyle Q_{yyi,j}=&\displaystyle q_{yyi,j}+\frac{\Delta x^2}{24}q_{yyxx}(x_{i},y_{j})+\frac{\Delta y^2}{24}q_{yyyy}(x_{i},y_{j})\\
&+\mathcal{O}(\Delta x^4+\Delta x^2\Delta y^2+\Delta y^4),\\
\displaystyle Q_{xxxxi,j}=&\displaystyle q_{xxxxi,j}+\frac{\Delta x^2}{24}q_{xxxxxx}(x_{i},y_{j})+\frac{\Delta y^2}{24}q_{xxxxyy}(x_{i},y_{j})\\
&+\mathcal{O}(\Delta x^4+\Delta x^2\Delta y^2+\Delta y^4),\\
\displaystyle Q_{yyyyi,j}=&\displaystyle q_{yyyyi,j}+\frac{\Delta x^2}{24}q_{yyyyxx}(x_{i},y_{j})+\frac{\Delta y^2}{24}q_{yyyyyy}(x_{i},y_{j})\\
&+\mathcal{O}(\Delta x^4+\Delta x^2\Delta y^2+\Delta y^4),\\
\displaystyle Q_{xxyyi,j}=&\displaystyle q_{xxyyi,j}+\frac{\Delta x^2}{24}q_{xxyyxx}(x_{i},y_{j})+\frac{\Delta y^2}{24}q_{xxyyyy}(x_{i},y_{j})\\
&+\mathcal{O}(\Delta x^4+\Delta x^2\Delta y^2+\Delta y^4).
\end{array}
\end{equation}
Therefore, if cell average  values $Q_{i,j}$ are available, we can approximate the second derivatives
$q_{xxi,j},q_{yyi,j}$, and the fourth derivatives $q_{xxxxi,j}$, $q_{xxyyi,j}$, $q_{yyyyi,j}$   to some order of accuracy  by using standard finite difference schemes for
 $Q_{xxi,j}, Q_{yyi,j}$, $Q_{xxxxi,j}$, $Q_{xxyyi,j}$, $Q_{yyyyi,j}$  in  (\ref{d2}). This is critical for the high-order conversion from average values to point values.

\subsubsection{Approximation of derivatives from point values}
\label{sec:2.2.1}
In order to get a six-order conversion formula from Eq.~(\ref{se3:eq5}) with points values of $q_{i,j}$ available, fourth-order accurate representations of $q_{xx}(x_{i},y_{j})$ and $q_{yy}(x_{i},y_{j})$ are required and can be obtained directly by using the standard fourth-order accurate finite difference scheme
\begin{equation}
\label{ptod4}
\begin{split}
\displaystyle q_{xx}(x_{i},y_{j})=\frac{1}{12\Delta x^2}\left(-q_{i-2,j}+16q_{i-1,j}-30q_{i,j}+16q_{i+1,j}-q_{i+2,j}\right)+\mathcal{O}(\Delta x^4),\\
\displaystyle q_{yy}(x_{i},y_{j})=\frac{1}{12\Delta y^2}\left(-q_{i,j-2}+16q_{i,j-1}-30q_{i,j}+16q_{i,j+1}-q_{i,j+2}\right)+\mathcal{O}(\Delta y^4).
\end{split}
\end{equation}
Only second-order accurate representations of $q_{xxxx}(x_{i},y_{j})$,  $ q_{xxyy}(x_{i},y_{j})$ and $q_{yyyy}(x_{i},y_{j})$ in Eq.~(\ref{se3:eq5}) are required  and can be obtained directly by
using the standard finite difference schemes
\begin{equation}
\label{ptod3}
\setlength{\arraycolsep}{0.5pt}
\begin{array}{ll}
\displaystyle q_{xxxx}(x_{i},y_{j})=&\displaystyle\frac{1}{\Delta x^4}\left(q_{i-2,j}-4q_{i-1,j}+6q_{i,j}-4q_{i+1,j}+q_{i+2,j}\right)+\mathcal{O}(\Delta x^2),\\
\displaystyle q_{yyyy}(x_{i},y_{j})=&\displaystyle\frac{1}{\Delta y^4}\left(q_{i,j-2}-4q_{i,j-1}+6q_{i,j}-4q_{i,j+1}+q_{i,j+2}\right)+\mathcal{O}(\Delta y^2),\\
\displaystyle q_{xxyy}(x_{i},y_{j})=&\displaystyle\frac{1}{\Delta x^2\Delta y^2}\left[(q_{i-1,j-1}+q_{i+1,j-1}-2q_{i,j-1})+(q_{i-1,j+1}+q_{i+1,j+1}\right.\\
&\displaystyle\left.-2q_{i,j+1})-2(q_{i-1,j}+q_{i+1,j}-2q_{i,j})\right]
+\mathcal{O}(\Delta x^2+\Delta y^2).
\end{array}
\end{equation}

\subsubsection{Approximation of derivatives from cell average values}
\label{sec:2.2.2}
In order to get six-order conversion formulas from  Eq.~(\ref{se3:eq5}) with  cell averaged values of $Q_{i,j}$ available,
we derive new fourth-order accurate approximations to $q_{xx}(x_{i},y_{j})$ and $q_{yy}(x_{i},y_{j})$ by replacing $q_{xxxxi,j}$, $q_{yyyyi,j}$ and $q_{xxyyi,j}$ in the first two equalities of (\ref{d2}) with $Q_{xxxx}(x_i,y_j)$, $Q_{yyyy}(x_i,y_j)$ and $Q_{xxyy}(x_i,y_j)$ in the last three equalities of (\ref{d2}), and then discretizing $Q_{xxi,j}$ and $Q_{yyi,j}$ with a fourth-order accurate finite difference  like (\ref{ptod4}) and $Q_{xxxx}(x_i,y_j)$, $Q_{yyyy}(x_i,y_j)$ and $Q_{xxyy}(x_i,y_j)$  with  second-order accurate finite differences like (\ref{ptod3}), giving
\begin{equation}
\label{atod2}
\setlength{\arraycolsep}{0.5pt}
\begin{array}{ll}
\displaystyle q_{xx}(x_{i},y_{j})&=\displaystyle Q_{xxi,j}-\frac{\Delta x^2}{24}Q_{xxxx}(x_{i},y_{j})-\frac{\Delta y^2}{24}Q_{xxyy}(x_{i},y_{j})
+\mathcal{O}(\Delta x^4+\Delta x^2\Delta y^2+\Delta y^4)\\
&\displaystyle =   \frac{1}{8\Delta x^2}\left(-Q_{i-2,j}+12Q_{i-1,j}-22Q_{i,j}+12Q_{i+1,j}-Q_{i+2,j}\right)\\
&\displaystyle~~~-\frac{1}{24\Delta x^2}\left[(Q_{i-1,j-1}+Q_{i+1,j-1}-2Q_{i,j-1})+(Q_{i-1,j+1}\right.\\
&\displaystyle~~~\left.+Q_{i+1,j+1}-2Q_{i,j+1})-2(Q_{i-1,j}+Q_{i+1,j}-2Q_{i,j})\right]\\
&\displaystyle~~~+\mathcal{O}(\Delta x^4+\Delta x^2\Delta y^2+\Delta y^4), \\
\displaystyle q_{yy}(x_{i},y_{j})&=\displaystyle Q_{yyi,j}-\frac{\Delta x^2}{24}Q_{yyxx}(x_{i},y_{j})-\frac{\Delta y^2}{24}Q_{yyyy}(x_{i},y_{j})
+\mathcal{O}(\Delta x^4+\Delta x^2\Delta y^2+\Delta y^4)\\
&\displaystyle = \frac{1}{8\Delta y^2}\left(-Q_{i,j-2}+12Q_{i,j-1}-22Q_{i,j}+12Q_{i,j+1}-Q_{i,j+2}\right)\\
&\displaystyle~~~-\frac{1}{24\Delta y^2}\left[(Q_{i-1,j-1}+Q_{i-1,j+1}-2Q_{i-1,j})+(Q_{i+1,j-1}\right.\\
&\displaystyle~~~\left.+Q_{i+1,j+1}-2Q_{i+1,j})-2(Q_{i,j-1}+Q_{i,j+1}-2Q_{i,j})\right]\\
&\displaystyle~~~+\mathcal{O}(\Delta x^4+\Delta x^2\Delta y^2+\Delta y^4).
\end{array}
\end{equation}

The second-order accurate representations of $q_{xxxx}(x_{i},y_{j})$, $q_{xxyy}(x_{i},y_{j})$  and $q_{yyyy}(x_{i},y_{j})$  in terms of cell average values can be obtained by using
the last three equalities in (\ref{d2}) and discretizing $Q_{xxxxi,j}$, $Q_{yyyyi,j}$ and $Q_{xxyyi,j}$ with the standard second-order accurate finite differences similar to (\ref{ptod3}),
\begin{equation}
\label{atod1}
\setlength{\arraycolsep}{0.5pt}
\begin{array}{ll}
\displaystyle q_{xxxx}(x_{i},y_{j})&=\displaystyle Q_{xxxxi,j}+\mathcal{O}(\Delta x^2+\Delta y^2)\\
&=\displaystyle\frac{1}{\Delta x^4}\left(Q_{i-2,j}-4Q_{i-1,j}+ 6Q_{i,j}-4Q_{i+1,j}+Q_{i+2,j}\right)+\mathcal{O}(\Delta x^2+\Delta y^2),\\
\displaystyle q_{yyyy}(x_{i},y_{j})&=\displaystyle Q_{yyyyi,j}+\mathcal{O}(\Delta x^2+\Delta y^2)\\
&=\displaystyle\frac{1}{\Delta y^4}\left(Q_{i,j-2}-4Q_{i,j-1}+6Q_{i,j}-4Q_{i,j+1}+Q_{i,j+2}\right)+\mathcal{O}(\Delta x^2+\Delta y^2),\\
\displaystyle q_{xxyy}(x_{i},y_{j})&=\displaystyle Q_{xxyyi,j}+\mathcal{O}(\Delta x^2+\Delta y^2)\\
&=\displaystyle\frac{1}{\Delta x^2\Delta y^2}\left[(Q_{i-1,j-1}+Q_{i+1,j-1}-2Q_{i,j-1})+(Q_{i-1,j+1}+Q_{i+1,j+1}\right.\\
&~~~\displaystyle\left.-2Q_{i,j+1})-2(Q_{i-1,j}+Q_{i+1,j}-2Q_{i,j})\right]+\mathrm{O}(\Delta x^2+\Delta y^2).\\
\end{array}
\end{equation}
We remark that one can get even higher order conversion  by the above process.

\subsection{Modified dimension-by-dimension FV WENO method}
\label{sec:2.3}
 The utilization of the conversion formulas between point values and cell average vales in the previous subsection suggests the following modified dimension-by-dimension FV WENO method on 3D Cartesian grids.
\paragraph{Algorithm: modified FV WENO method with sixth-order conversion formulas}
\begin{enumerate}[(1)]
\item [1.] Compute face-averaged values of the conserved quantities at cell interfaces using the one-dimensional WENO reconstruction, i.e. compute
$$U^{\pm}_{i+\frac{1}{2},j,k}(t), \quad U^{\pm}_{i,j+\frac{1}{2},k}(t), \quad  U^{\pm}_{i,j,k+\frac{1}{2}}(t) $$
 at all grid cell interfaces including several layers of ghost cell faces outside the computational domain. Notice that boundary conditions are applied to obtain the conserved quantities $U_{i,j,k}$ in ghost cells so that we can implement the 1D WENO reconstruction for the ghost cell interfaces.
\item [2.] Compute point values of the conserved variables at the centers of cell interfaces, i.e. compute
$$u^{\pm}_{i+\frac{1}{2},j,k}(t), \quad u^{\pm}_{i,j+\frac{1}{2},k}(t), \quad  u^{\pm}_{i,j,k+\frac{1}{2}}(t) $$
at all interface centers including several layers of ghost interface centers using the conversion formula (\ref{se3:eq5}) as substantialized by Eq.~(\ref{u6}).
\item [3.] Compute point values of the numerical fluxes at all interface centers including ghost ones, i.e.
\begin{equation*}
\begin{split}
&\hat{f}_{i+\frac{1}{2},j,k}(t)=\mathcal{F}\left(u^-_{i+\frac12,j,k}, u^+_{i+\frac12,j,k}\right), \quad\hat{g}_{i,j+\frac{1}{2},k}(t)=\mathcal{F}\left(u^-_{i,j+\frac12,k}, u^+_{i,j+\frac12,k}\right),  \\
&\hat{h}_{i,j,k+\frac{1}{2}}(t)=\mathcal{F}\left(u^-_{i,j,k+\frac12}, u^+_{i,j,k+\frac12}\right)
\end{split}
\end{equation*}
using a numerical flux function like the Lax-Friedrichs flux or HLLC flux.
\item [4.]  Compute averaged values of the numerical flux at  grid cell
interfaces
$$\hat{F}_{i+\frac{1}{2},j,k}(t),\quad \hat{G}_{i,j+\frac{1}{2},k}(t), \quad \hat{H}_{i,j,k+\frac{1}{2}}$$
using the conversion formula (\ref{se3:eq5}) as substantialized by Eq.~(\ref{f6}).
\item [5.] Solve the semi-discrete system  (\ref{se2:eq3}) using a high-order accurate Runge-Kutta method.
\end{enumerate}

In this paper, we consider the standard dimension-by-dimension WENO method and the modified FV WENO method with sixth-order conversion formulas.
For the Euler equations, in step 1 we use the characteristic variables  for the WENO reconstruction \cite{shu96}.
\paragraph{Classical method:}  The standard dimension-by-dimension FV WENO method.
\paragraph{Modified method:}
Point values of the conserved variables at the center points
of cell interfaces in  step 2 are computed by using the following sixth-order accurate average-to-point conversion formula, which is obtained by substituting (\ref{atod2}) and (\ref{atod1}) into  (\ref{se3:eq5}),
\begin{equation}
\label{u6}
\setlength{\arraycolsep}{0.5pt}
\begin{array}{lll}
\displaystyle u^{\pm}_{i+\frac{1}{2},j,k}&=&\displaystyle U^{\pm}_{i+\frac{1}{2},j,k}-\frac{1}{1920} \left( -9U^{\pm}_{i+\frac{1}{2},j-2,k} +116U^{\pm}_{i+\frac{1}{2},j-1,k}-214U^{\pm}_{i+\frac{1}{2},j,k} \right.\\
& &\displaystyle\left.+116U^{\pm}_{i+\frac{1}{2},j+1,k}-9U^{\pm}_{i+\frac{1}{2},j+2,k}\right)
-\frac{1}{1920}\left(-9U^{\pm}_{i+\frac{1}{2},j,k-2}\right.\\
& &\displaystyle\left.+116U^{\pm}_{i+\frac{1}{2},j,k-1}- 214U^{\pm}_{i+\frac{1}{2},j,k} +116U^{\pm}_{i+\frac{1}{2},j,k+1}-9U^{\pm}_{i+\frac{1}{2},j,k+2}\right)\\
& &\displaystyle+\frac{1}{576}\left[ \left(U^{\pm}_{i+\frac{1}{2},j-1,k-1}+U^{\pm}_{i+\frac{1}{2},j+1,k-1}-2U^{\pm}_{i+\frac{1}{2},j,k-1} \right) \right.\\
& &\displaystyle+\left(U^{\pm}_{i+\frac{1}{2},j-1,k+1}+U^{\pm}_{i+\frac{1}{2},j+1,k+1}-2U^{\pm}_{i+\frac{1}{2},j,k+1} \right) \\
& &\displaystyle\left.-2\left(U^{\pm}_{i+\frac{1}{2},j-1,k}+U^{\pm}_{i+\frac{1}{2},j+1,k}-2U^{\pm}_{i+\frac{1}{2},j,k} \right) \right],
\end{array}
\end{equation}
and analogously for  \(u^{\pm}_{i,j+\frac{1}{2},k}(t)$ and $ u^{\pm}_{i,j,k+\frac{1}{2}}(t)\).

Averaged values of the numerical flux at cell
interfaces in step 4 are computed by using the following sixth-order accurate point-to-average conversion formula obtained by substituting (\ref{ptod4}) and  (\ref{ptod3}) into (\ref{se3:eq5}),
\begin{equation}
\label{f6}
\setlength{\arraycolsep}{0.5pt}
\begin{array}{lll}
\displaystyle \hat{F}_{i+\frac{1}{2},j,k}&=&\displaystyle \hat{f}_{i+\frac{1}{2},j,k}+\frac{1}{5760} \left(-17\hat{f}_{i+\frac{1}{2},j-2,k}+ 308\hat{f}_{i+\frac{1}{2},j-1,k}-582\hat{f}_{i+\frac{1}{2},j,k}\right.\\
& &\displaystyle\left. +308\hat{f}_{i+\frac{1}{2},j+1,k}-17\hat{f}_{i+\frac{1}{2},j+2,k}\right)
+\frac{1}{5760}\left(-17\hat{f}_{i+\frac{1}{2},j,k-2}\right.\\
& &\displaystyle\left.+308\hat{f}_{i+\frac{1}{2},j,k-1}-582\hat{f}_{i+\frac{1}{2},j,k} +308\hat{f}_{i+\frac{1}{2},j,k+1}-17\hat{f}_{i+\frac{1}{2},j,k+2}\right)\\
& &\displaystyle+\frac{1}{576}\left [ \left(\hat{f}_{i+\frac{1}{2},j-1,k-1}+\hat{f}_{i+\frac{1}{2},j+1,k-1}-2\hat{f}_{i+\frac{1}{2},j,k-1} \right) \right.\\
& &\displaystyle+\left(\hat{f}_{i+\frac{1}{2},j-1,k+1}+\hat{f}_{i+\frac{1}{2},j+1,k+1}-2\hat{f}_{i+\frac{1}{2},j,k+1} \right) \\
& &\displaystyle\left.-2\left(\hat{f}_{i+\frac{1}{2},j-1,k}+\hat{f}_{i+\frac{1}{2},j+1,k}-2\hat{f}_{i+\frac{1}{2},j,k} \right) \right],
\end{array}
\end{equation}
and analogously for $\hat{G}_{i,j+\frac{1}{2},k}$ and $ \hat{H}_{i,j,k+\frac{1}{2}}$.

In Table \ref{table:t1}, we summarize the expected convergence rates of two different methods
for the approximation of linear and nonlinear smooth problems.

\begin{table}
\centering
\caption{Predicted convergence rate for
 different methods}
\label{table:t1}
\begin{tabular}{l  c  c  c  c  c}
\hline
\multirow{2}{*}{Method} & \multicolumn{2}{c}{WENO5+RK5}&  & \multicolumn{2}{c}{WENO7+RK7}\\
\cline{2-3}
\cline{5-6}
\multirow{2}{*}{} & Linear  & Nonlinear & &  Linear & Nonlinear\\
\midrule[1pt]
Classical method   & 5        &   2       & &    7    &    2     \\
Modified method    & 5        &   5       & &    6    &    6     \\
\hline
\end{tabular}
\end{table}

\begin{figure}[hp]
  \begin{minipage}[t]{0.5\linewidth}
    \centering
    \includegraphics[width=2.5in,height=1.8in, bb=0 0 671 395]{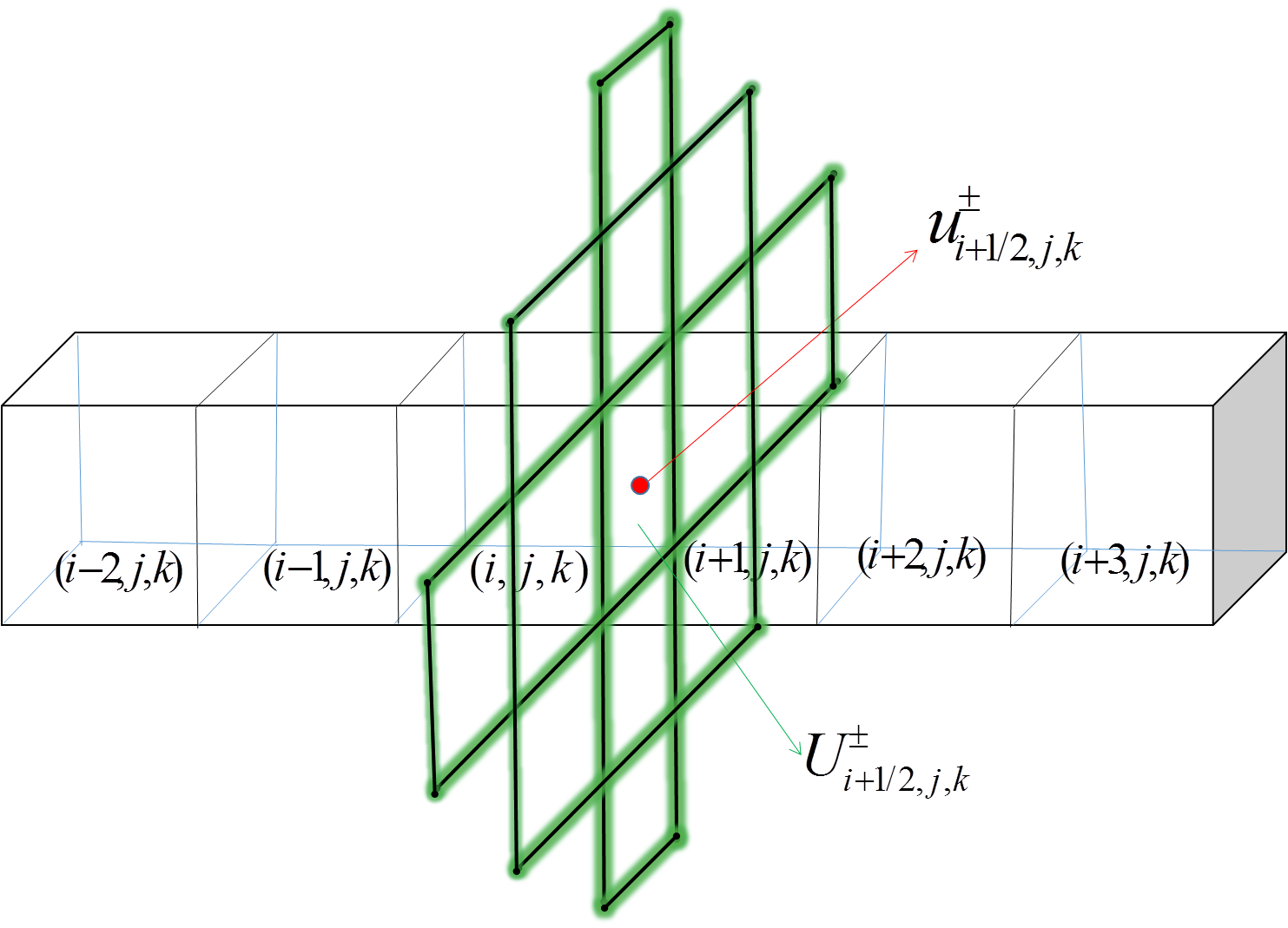}
    \caption{The stencil for computing the face center point values  $u^{\pm}_{i+\frac{1}{2},j,k}$ marked as the red point. The green face cells denote face averaged values $U_{i+1/2,j,k}$ obtained from 1D WENO reconstructions and used to compute  $u^{\pm}_{i+\frac{1}{2},j,k}$ by Eq.~(\ref{u6}). Only one WENO stencil along the $i$ direction is shown for clarity.}
    \label{fig:fig21}
  \end{minipage}
  \makebox{}
  \begin{minipage}[t]{0.5\linewidth}
    \centering
    \includegraphics[width=2.2in,height=1.8in, bb=0 0 586 474]{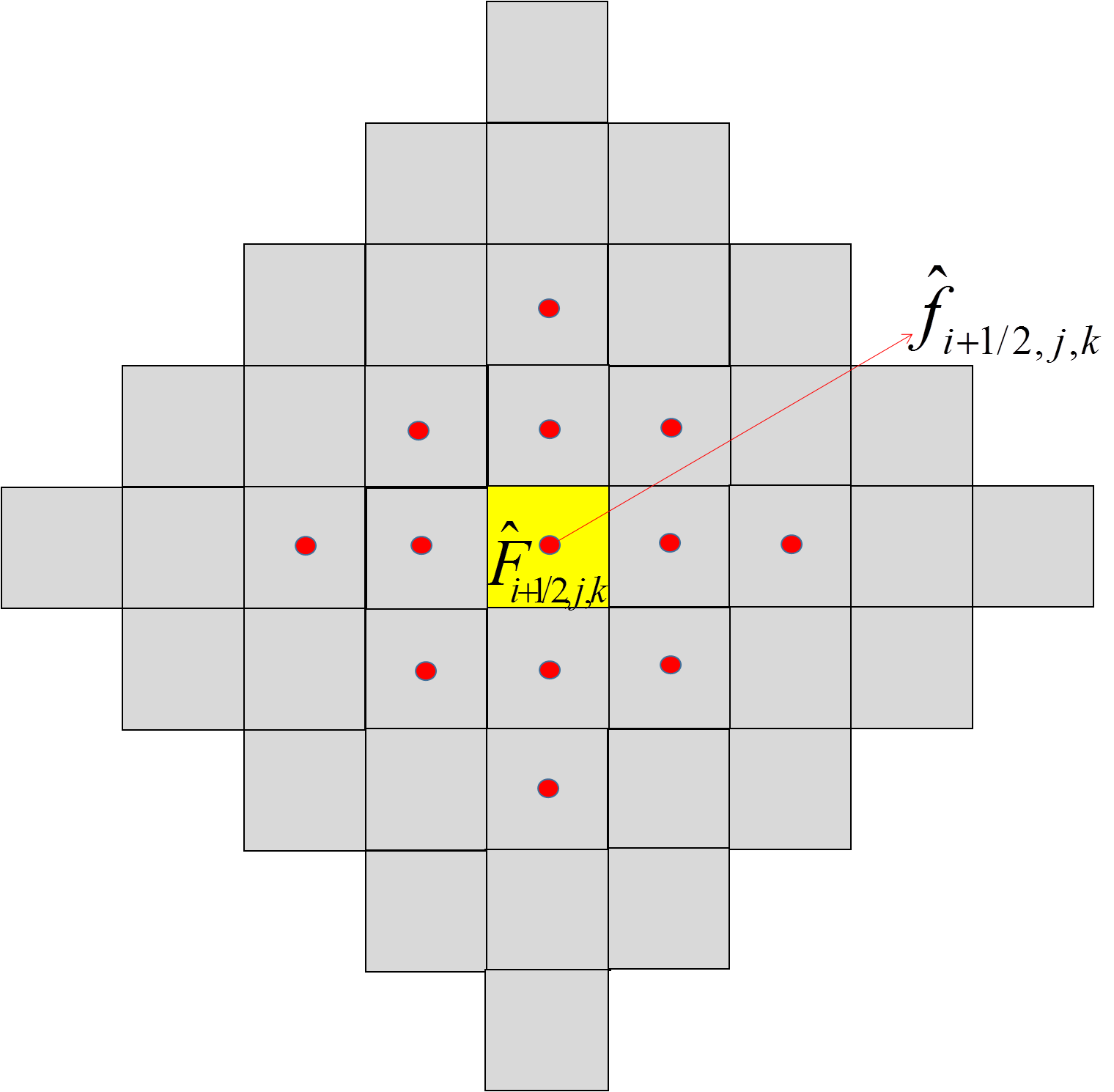}
    \caption{The $yz$ plane view of the stencil for computing the face averaged value $\hat{F}_{i+\frac{1}{2},j,k}$ marked as yellow cell.  The red points denote the point values $\hat{f}_{i+1/2,j,k}$ used to compute $\hat{F}_{i+\frac{1}{2},j,k}$  by Eq.~(\ref{f6}). The shaded cells represent the involved face averaged values $U_{i+1/2,j,k}$ used to compute the red point values $u^{\pm}_{i+\frac{1}{2},j,k}$ thus $\hat{f}_{i+1/2,j,k}$.}
   \label{fig:fig22}
  \end{minipage}
\end{figure}

For the modified method, Fig.~\ref{fig:fig21} shows the stencil for computing the face center point values of the conserved quantities $u^\pm_{i+1/2,j,k}$ using (\ref{u6}). This transformation requires several face averaged values  $U_{i+1/2,j,k}$ of the conserved quantities on the $i+1/2$ plane marked as green face cells, each of which is obtained from the 1D WENO reconstruction in the $i$ direction using cell averages from $i-2$ to $i+3$ as shown for the WENO5 case.
Fig.~\ref{fig:fig22}  shows the stencil for computing the face averaged value of the flux $\hat{F}_{i+1/2,j,k}$ using (\ref{f6}). This transformation requires several point values of fluxes $\hat{f}_{i+1/2,j,k}$ marked as red points, each of which comes from $u^\pm_{i+1/2,j,k}$ that further involves several $i$-direction WENO stencils required by Eq.~(\ref{u6}). The projections of all the involved $i$-direction WENO stencils onto the $i+1/2$ plane required for computing the
face averaged flux $\hat{F}_{i+1/2,j,k}$  are marked as shaded cells.
\newpage

\section{ Numerical results}
\label{sec:3}
In this section, several 3D numerical examples are used to compare the performance of the present modified FV WENO method and the classic  method.

In the following tables to show convergence studies, the $\|\cdot\|_1$ norm of the error denotes the quantity $\sum\limits_{i,j,k} |\bar{u}_{i,j,k}-\bar{u}^{\text{exact}}_{i,j,k} | \times \Delta x \Delta y \Delta z$, where $\bar{u}$ represents the cell average value. We compute the experimental order
of convergence (EOC) using the formula
\begin{equation}
\displaystyle\text{EOC}=\frac{\displaystyle  \log\left(\|\bar{u}_{m}- \bar{u}_{\text{exact}}\|_{1}/\|\bar{u}_{2m}- \bar{u}_{\text{exact}}\|_{1}  \right)}{\log 2},
\end{equation}
where the index $m$ indicates the number of grid cells in the $x, y$ and $z$ direction. We use fifth or seventh order accuracy WENO-Z reconstruction \cite{wenoz,wenoz7,highweno} with the parameters $q=2$ and $\epsilon=10^{-14}$. For the temporal discretization, in order to match the order of spatial accuracy, fifth or seventh explicit Runge-Kutta  schemes (see Appendix 1 of Ref.~\cite{improved14}) are used.   In all computations, the time steps used correspond to $\text{CFL} =0.5$.

\subsection{Scalar hyperbolic problems}
\label{sec:3.1}
\subsubsection{3D linear advection equation}
\label{sec:3.1.1}
We consider the 3D linear advection problem \cite{zhoutie01} given by
\begin{equation}
\label{3dlinear}\left\{
\begin{array}{l}
u_t+u_x+u_y+u_z=0,\\
u(x,y,z,0)=\displaystyle \sin\left(\frac{\pi}{2}(x+y+z)\right),
\end{array}\right. -2\leq x,y,z\leq 2,
\end{equation}
with  periodic boundary conditions.  The exact solution is $u(x,y,z,t)=\sin\left(\frac{\pi}{2}(x+y\right.$ $\left.+z-3t)\right)$. The final time is $T=1.0$. The Lax-Friedrichs flux is used.

In Table \ref{table:t2}, we show the $\|\cdot\|_1$ norm of the error  and orders of grid convergence for the problem  (\ref{3dlinear}) by using the fifth order WENO-Z reconstruction and RK5 time stepping scheme for the two FV WENO methods. As expected for the linear problem,  the classical method converges with nearly fifth order accuracy, The full fifth order of convergence of the WENO-Z
reconstruction is also nearly retained by modified method as expected in Table \ref{table:t1}. We can also observe that the absolute errors for the classical method and modified method are nearly the same.

In Table \ref{table:t3}, we show the results computed by using the two methods with the seventh order WENO-Z reconstruction and RK7 time stepping scheme for the same problem  (\ref{3dlinear}). Again we observe that the order of convergence is as expected in Table \ref{table:t1} for each method.
\begin{table}
\centering
\caption{Convergence study for problem (\ref{3dlinear}) with 5th order WENO-Z reconstruction and RK5.}
\label{table:t2}
\begin{tabular}{l  c  c  c  c  c}
 \hline
 \multirow{2}{*}{Grid} & \multicolumn{2}{c}{Classical method}&  & \multicolumn{2}{c}{Modified method} \\
 \cline{2-3}
 \cline{5-6}
\multirow{2}{*}{} & $\|\bar{u}-\bar{u}_{\text{exact}}\|_{1}$  & EOC &  & $\|\bar{u}-\bar{u}_{\text{exact}}\|_{1}$ &  EOC \\
\midrule[1pt]
$10^{3}$ & $4.6824\times 10^{-3}$ &      & & $4.7219\times 10^{-3}$ &  \\
$20^{3}$ & $1.5639\times 10^{-4}$ & 4.90 & & $1.5641\times 10^{-4}$ & 4.91\\
$40^{3}$ & $5.0886\times 10^{-6}$ & 4.94 & & $5.0886\times 10^{-6}$ & 4.94\\
$80^{3}$ & $1.6068\times 10^{-7}$ & 4.98 & & $1.6068\times 10^{-7}$ & 4.98\\
\hline
\end{tabular}
\end{table}

\begin{table}
\centering
\caption{Convergence study for problem (\ref{3dlinear}) with 7th order WENO-Z reconstruction and RK7.}
\label{table:t3}
\begin{tabular}{l  c  c   c  c  c}
 \hline
 \multirow{2}{*}{Grid} & \multicolumn{2}{c}{Classical method}&  & \multicolumn{2}{c}{Modified method} \\
 \cline{2-3}
 \cline{5-6}
\multirow{2}{*}{} & $\|\bar{u}-\bar{u}_{\text{exact}}\|_{1}$  & EOC &  & $\|\bar{u}-\bar{u}_{\text{exact}}\|_{1}$ &  EOC \\
\midrule[1pt]
$10^{3}$ & $4.7178\times 10^{-4}$ &      & & $6.4394\times 10^{-4}$ &  \\
$20^{3}$ & $3.1288\times 10^{-6}$ & 7.24 & & $5.0894\times 10^{-6}$ & 6.98\\
$40^{3}$ & $2.5216\times 10^{-8}$ & 6.96 & & $6.3723\times 10^{-8}$ & 6.32\\
$80^{3}$ & $1.9853\times 10^{-10}$& 6.99 & & $9.2312\times 10^{-10}$& 6.11\\
\hline
\end{tabular}
\end{table}

\subsubsection{3D Burgers equation}
\label{sec:3.1.2}
We consider the 3D Burgers equation problem \cite{toro04,zhoutie01}
\begin{equation}\left\{
\begin{split}
\displaystyle u_t+\left(\frac{u^{2}}{2}\right)_x+\left(\frac{u^{2}}{2}\right)_y+ \left(\frac{u^{2}}{2}\right)_z=0,\\
u(x,y,z,0)=0.5+\sin\left(\frac{\pi}{3}(x+y+z)\right),
\end{split}
\right.  ~~ -3\leq x,y,z\leq 3.
\label{3dburger}
\end{equation}
In this test, the boundary conditions are periodic, and the Lax-Friedrichs flux is used. The final time is set to $T=0.1$ before the discontinuity occurs. Here, the analytic solution is  $u(x,y,z,t)=0.5+\sin\left(\frac{\pi}{3}(x+y+z-3ut)\right)$. The cell average values of the exact solution are computed by using the Gauss quadrature formula with ninth order accuracy.

In Tables \ref{table:t4} and \ref{table:t5}, we show the error and the grid convergence rates for problem (\ref{3dburger}) by using the two methods.
The fifth order WENO-Z reconstruction together with RK5 is used for Table \ref{table:t4} and the seventh order WENO-Z reconstruction with RK7 is used for Table \ref{table:t5}. As expected, the convergence rate for the classical method is only second order for nonlinear problems. In comparison, the orders of convergence for the modified method could reach fifth in Table \ref{table:t4} and sixth in Table \ref{table:t5} respectively on refined grids.
\begin{table}
\centering
\caption{Convergence study for problem (\ref{3dburger}) with 5th order WENO-Z reconstruction and RK5.}
\label{table:t4}
\begin{tabular}{l  c  c  c  c  c}
 \hline
 \multirow{2}{*}{Grid} & \multicolumn{2}{c}{Classical method}&  & \multicolumn{2}{c}{Modified method} \\
 \cline{2-3}
 \cline{5-6}
\multirow{2}{*}{} & $\|\bar{u}-\bar{u}_{\text{exact}}\|_{1}$  & EOC & & $\|\bar{u}-\bar{u}_{\text{exact}}\|_{1}$ &  EOC \\
\midrule[1pt]
$10^{3}$ & $5.8528\times 10^{-3}$ &      & & $6.5145\times 10^{-4}$ &  \\
$20^{3}$ & $1.6359\times 10^{-3}$ & 1.84 & & $7.0417\times 10^{-5}$ & 3.21\\
$40^{3}$ & $4.5129\times 10^{-4}$ & 1.86 & & $3.1143\times 10^{-6}$ & 4.50\\
$80^{3}$ & $1.1381\times 10^{-4}$ & 1.99 & & $1.0734\times 10^{-7}$ & 4.86\\
$160^{3}$& $2.8503\times 10^{-5}$ & 2.00 & & $3.4787\times 10^{-9}$ & 4.95\\
\hline
\end{tabular}
\end{table}
\begin{table}
\centering
\caption{Convergence study for problem (\ref{3dburger}) with 7th order WENO-Z reconstruction and RK7.}
\label{table:t5}
\begin{tabular}{l  c  c  c  c  c}
 \hline
 \multirow{2}{*}{Grid} & \multicolumn{2}{c}{Classical method}&  & \multicolumn{2}{c}{Modified method} \\
 \cline{2-3}
 \cline{5-6}
\multirow{2}{*}{} & $\|\bar{u}-\bar{u}_{\text{exact}}\|_{1}$  & EOC & &$\|\bar{u}-\bar{u}_{\text{exact}}\|_{1}$ &  EOC \\
\midrule[1pt]
$10^{3}$ & $5.7287\times 10^{-3}$ &      & & $4.8067\times 10^{-4}$ &  \\
$20^{3}$ & $1.6668\times 10^{-3}$ & 1.78 & & $1.3965\times 10^{-5}$ & 5.11\\
$40^{3}$ & $4.5324\times 10^{-4}$ & 1.89 & & $3.9716\times 10^{-7}$ & 5.14\\
$80^{3}$ & $1.1389\times 10^{-4}$ & 1.99 & & $3.8968\times 10^{-9}$ & 6.67\\
$160^{3}$& $2.8505\times 10^{-5}$ & 2.00 & & $5.3592\times 10^{-11}$& 6.05\\
\hline
\end{tabular}
\end{table}

\subsection{3D Euler equations}
\label{sec:3.2}
 In this subsection we use the 3D Euler equations of gas dynamics
\begin{equation}\partial_{t}\left(\begin{matrix} \rho\\ \rho u\\ \rho v \\ \rho w \\ E \end{matrix}\right)+ \partial_{x}\left(\begin{matrix} \rho u\\ \rho u^{2}+p\\ \rho uv \\ \rho uw \\ u(E+p)\end{matrix}\right)+ \partial_{y}\left(\begin{matrix} \rho v\\ \rho uv\\ \rho v^{2}+p \\ \rho vw \\ v(E+p) \end{matrix}\right)+ \partial_{z}\left(\begin{matrix} \rho w\\ \rho uw\\ \rho vw \\ \rho w^{2}+p \\ w(E+p) \end{matrix}\right)=0
\label{3deuler}
\end{equation}
as our model problem with the ideal gas equation of state
 $$E=\frac{p}{\gamma-1}+\frac{1}{2} \rho(u^{2}+v^{2}+w^{2}). $$
The initial values and boundary conditions will be specified below for each test problem. We always set  $\gamma=1.4$, and
use the characteristic variables for the WENO reconstruction and the  HLLC flux for the numerical flux.
\subsubsection{Linear problem}
\label{sec:3.2.1}
We consider the periodic solutions \cite{shu11} of the Euler equations (\ref{3deuler}).  The initial values are given by
\begin{equation}\left\{
\begin{split}
\rho(x,y,z,0)&=\displaystyle 1+0.2\sin\left(\frac{\pi}{3}(x+y+z)\right),\\
p(x,y,z,0)&=1,\\
u(x,y,z,0)&=v(x,y,z,0)=w(x,y,z,0)=1,
\end{split}\right. ~~-3\leq x,y,z\leq 3.
\label{3deularlinear}
\end{equation}
 Periodic boundary conditions are applied in this test. The exact solution of density is $\rho(x,y,z,t)$ $=1+0.2\sin\left(\frac{\pi}{3}(x+y+z-3t)\right)$. The final time is $T=1.0$.

\begin{table}
\centering
\caption{Convergence study for problem (\ref{3deularlinear}) with 5th order WENO-Z reconstruction and RK5.}
\label{table:t6}
\begin{tabular}{l  c  c  c  c  c}
 \hline
 \multirow{2}{*}{Grid} & \multicolumn{2}{c}{Classical method}& & \multicolumn{2}{c}{Modified method} \\
 \cline{2-3}
 \cline{5-6}
\multirow{2}{*}{} & $\|\bar{\rho}-\bar{\rho}_{\text{exact}}\|_{1}$  & EOC &  &  $\|\bar{\rho}-\bar{\rho}_{\text{exact}}\|_{1}$ &  EOC \\
\midrule[1pt]
$10^{3}$ & $1.4454\times 10^{-3}$ &      & & $1.5166\times 10^{-3}$ &  \\
$20^{3}$ & $4.5141\times 10^{-5}$ & 5.00 & & $4.8440\times 10^{-5}$ & 4.97\\
$40^{3}$ & $1.4121\times 10^{-6}$ & 5.00 & & $1.5275\times 10^{-6}$ & 4.99\\
$80^{3}$ & $4.3960\times 10^{-8}$ & 5.00 & & $4.7722\times 10^{-8}$ & 5.00\\
\hline
\end{tabular}
\end{table}
\begin{table}
\centering
\caption{Convergence study for problem (\ref{3deularlinear}) with 7th order WENO-Z reconstruction and RK7.}
\label{table:t7}
\begin{tabular}{l  c  c  c  c  c}
 \hline
 \multirow{2}{*}{Grid} & \multicolumn{2}{c}{Classical method}&  & \multicolumn{2}{c}{Modified method} \\
 \cline{2-3}
 \cline{5-6}
\multirow{2}{*}{} & $\|\bar{\rho}-\bar{\rho}_{\text{exact}}\|_{1}$  & EOC & & $\|\bar{\rho}-\bar{\rho}_{\text{exact}}\|_{1}$ &  EOC \\
\midrule[1pt]
$10^{3}$ & $1.3190\times 10^{-4}$ &      & & $1.6761\times 10^{-4}$ &  \\
$20^{3}$ & $9.4877\times 10^{-7}$ & 7.12 & & $1.2019\times 10^{-6}$ & 7.12\\
$40^{3}$ & $7.4890\times 10^{-9}$ & 6.99 & & $1.1973\times 10^{-8}$ & 6.65\\
$80^{3}$ & $6.0914\times 10^{-11}$& 6.94 & & $1.4822\times 10^{-10}$& 6.34 \\
\hline
\end{tabular}
\end{table}

Tables \ref{table:t6} and \ref{table:t7} show the errors and numerical orders of accuracy of the density for problem (\ref{3deularlinear}) computed by using the two different FV WENO methods with the 5th order WENO-Z and 7th  WENO-Z reconstruction, respectively. In Table \ref{table:t6}, the errors of the classical method and modified method are very close. The numerical orders of accuracy are both five with 5th order WENO-Z reconstruction, which verify the theoretical prediction.
In Table \ref{table:t7}, the two FV WENO methods with the 7th WENO-Z reconstruction show similar trends like Table \ref{table:t6}. The classical method gives nearly seventh order of accuracy. And the modified method has well attained the theoretical sixth order of accuracy.

\subsubsection{Spherical Riemann problem}
\label{sec:3.2.2}
We test our modified FV WENO method for a three-dimensional spherical Riemann problem \cite{bump1,tu2007} between two parallel walls at $z=0$ and $z=1$ to observe whether the proposed method can work well for problems with discontinuities. The sphere is centered at $(0,0,0.4)$ with radius $r=0.2$. Initially the gas is at rest with density $\rho=1.0$ everywhere and pressure
\begin{equation*}
\left\{\begin{array}{ll} p=1 & \text{if~~}  r>0.2\\
p=5 & \text{else}   \end{array}\right.
\end{equation*}

The evolution of the flow field will remain cylindrically symmetric, thus a quarter computational domain is chosen to be $(x,y,z)\in[0,1.5]\times[0,1.5]\times[0,1]$. The grid is $150\times 150\times 100$. Reflective boundary conditions are imposed on the walls $z=0$ and $z=1$ and symmetric boundary conditions are used on symmetric planes $x=0$ and $y=0$. The other boundaries are taken to be outflow conditions.

\begin{figure}[hp!]
   \centering
  \begin{minipage}[t]{1\linewidth}
    \centering
    \includegraphics[height =1.8 in, bb=0 0 2000 869]{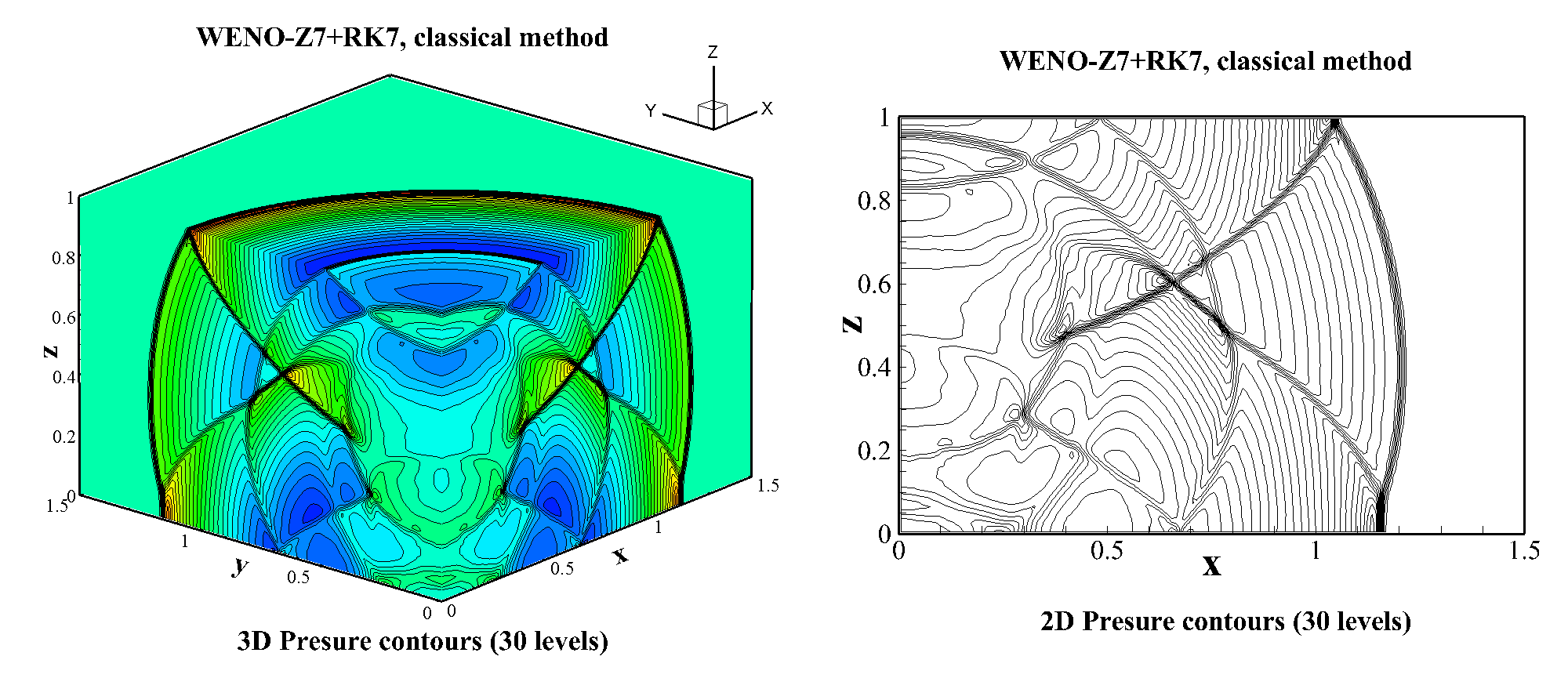}
  \end{minipage}
  \begin{minipage}[t]{1\linewidth}
    \centering
    \includegraphics[height =1.8 in,  bb=0 0 2000 869]{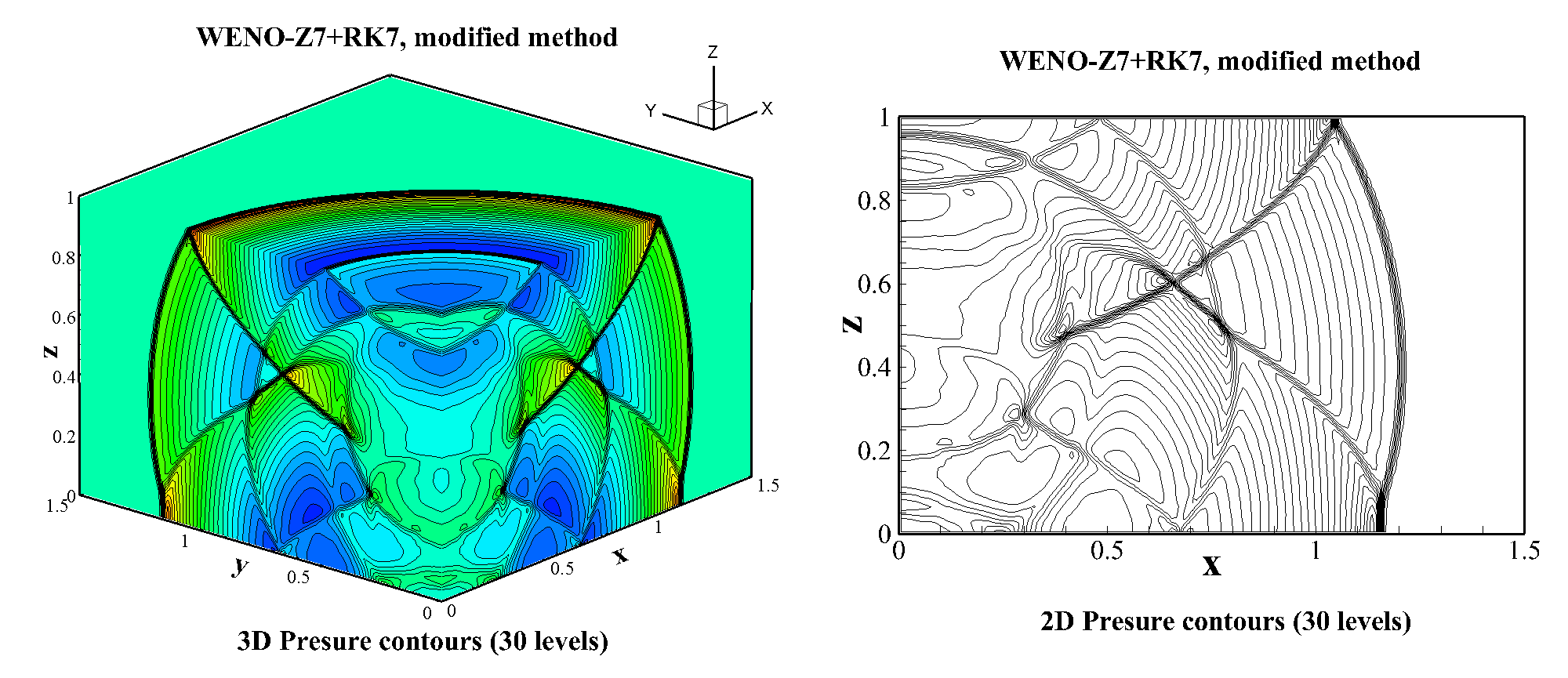}
  \end{minipage}
  \caption{Comparison of the classic  method with the modified method for the spherical Riemann problem at $T=0.7$ on the grid of $150\times150\times100$. 30 equally pressure contours from 0.781 to 1.475. Left: 3D view; Right: the $x-z$ plane view. Here WENO-Z7 with RK7 are used.
 }
\label{fig:f1}
\end{figure}

Fig.~\ref{fig:f1} shows the results for the spherical Riemann problem at time $T=0.7$ computed by using the two different methods. We can observe that the main features of the solution are the interactions between a strong outward moving shock wave, an outgoing imploded shock wave, and the walls. 
All the results computed by the two methods are in good agreement with other simulations \cite{bump1,tu2007}. We also observe that the modified method gives as good results as the classical method, which verifies that the central difference-based conversion formulas do not cause any numerical difficulty for problems with strong discontinuities. This was also noted in Ref.~\cite{improved14}.

\begin{table}[htp!]
\centering
\caption{Comparison for the average time consumption of one iteration between the classical and the modified methods with WENO-Z5+RK5 for the spherical Riemann problem. The ratio is relative to the classical method.}
\label{table:t8}
\begin{tabular}{l  l  c  c   l  c}
 \hline
 \multirow{2}{*}{Grid} & \multicolumn{2}{c}{Classical method}&   & \multicolumn{2}{c}{Modified method} \\
 \cline{2-3}
 \cline{5-6}
\multirow{2}{*}{} & Time (s)  & Ratio &  &   Time (s) &   Ratio \\
\midrule[1pt]
$37\times37\times25$   &0.4089  & 1.00  & & 0.8269   &2.03 \\
$75\times75\times50$   &3.1241  & 1.00  & & 6.1852   &1.98 \\
$150\times150\times100$&22.6309 & 1.00  & & 43.5660  &1.93 \\
\hline
\multicolumn{2}{l}{Average ratio}  & 1.00 & &        &1.98 \\
\hline
\end{tabular}
\end{table}
\begin{table}[htp!]
\centering
\caption{Comparison for the average time consumption of one iteration between the classical and the modified methods with WENO-Z7+RK7 for the spherical Riemann problem. The ratio is relative to the classical method.}
\label{table:t9}
\begin{tabular}{l  l  c  c  l  c}
 \hline
 \multirow{2}{*}{Grid} & \multicolumn{2}{c}{Classical method}& & \multicolumn{2}{c}{Modified method}\\
 \cline{2-3}
 \cline{5-6}
\multirow{2}{*}{} & Time (s)  & Ratio & &   Time (s) &  Ratio \\
\midrule[1pt]
$37\times37\times25$   &0.8548  & 1.00 & & 1.4717   &1.72 \\
$75\times75\times50$   &7.1570  & 1.00 & & 10.3105  &1.44 \\
$150\times150\times100$&54.2467 & 1.00 & & 77.8363  &1.43 \\
\hline
\multicolumn{2}{l}{Average ratio}     & 1.00 & &    &1.53 \\
\hline
\end{tabular}
\end{table}

Finally, in Tables \ref{table:t8} and \ref{table:t9} we show the timing results of the two FV WENO methods combined with the fifth or seventh order accuracy WENO-Z reconstruction for calculating the spherical Riemann problem. The codes have been parallelized using OpenMP and run by using 8 threads on a machine with 12 Intel Xeon(R) X5675 3.07~GHz CPU cores.  We see that the computational costs of the modified method increase on average by 98\% for the fifth order WENO-Z reconstruction, 53\% for the seventh order WENO-Z reconstruction. Moreover, the ratio is smaller for the more expensive seventh order WENO-Z7 reconstruction. This is because the additional computations required by the modified method are independent of the chosen reconstruction, which agrees well with the timing results \cite{improved14}.

\section{ Conclusions}
\label{sec:4}
We have further developed the modified dimension-by-dimension finite volume WENO method on three-dimensional Cartesian grids for nonlinear hyperbolic conservation laws. Specifically, we present  a sixth-order accurate formula for conversion from face average values to point values of the conserved quantities, and the underlying derivation process can be extended to even higher order of accuracy.  For three-dimensional problems, the computational cost of the modified FV WENO method is shown to be 1.25-2.03 times of that of the standard dimension-by-dimension finite volume  WENO method,  while other three-dimensional finite volume WENO methods on Cartesian grids \cite{yuan2011,toro04,hushu02} is   7--10 times of that as shown in Refs.~\cite{toro04,zhoutie01,dengxiaogang17}. The numerical tests show that the modified FV-WENO method, unlike the standard dimension-by-dimension finite volume  WENO method,
 retains the full spatial order of accuracy  when applied to smooth three-dimensional  nonlinear problems and is efficient and robust for calculating non-smooth nonlinear problems with shocks.
\section*{Acknowledgments}
This work is supported by Natural Science Foundation of China ($11321061$, $11261160486$, and 91641107), Fundamental Research of Civil Aircraft (MJ-F-2012-04).



\begin{thebibliography}{99}
\bibitem{shu96}G. S. Jiang and C. W. Shu, Efficient implementation of weighted ENO schemes,
J. Comput. Phys.  126 (1996), 202-228.
\bibitem{shu97}C. W. Shu, Essentially non-oscillatory and weighted essentially
non-oscillatory schemes for hyperbolic conservation laws, NASA/CR-97-206253, ICASE Report NO. 97-65, 1997.
\bibitem{shu11} R. Zhang, M. P. Zhang and C. W. Shu, On the order of accuracy and numerical performance of two classes of finite volume WENO schemes, Commun. Comput. Phys.  9 (2011), 807-827.
\bibitem{improved14}P. Buchm{\"{u}}ller and C. Helzel, Improved accuracy of high-order WENO finite volume methods on Cartesian grids,
J. Sci. Comput.  61 (2014), 343-368.

\bibitem{barth90}
T. J. Barth, P. Frederickson, High order solution of the Euler equations on unstructured grids using quadratic reconstruction. In: Proceedings of the twenty eighth aerospace science meeting, AIAA Paper No. 90-0013, 1990.

\bibitem{hushu99}C. Q. Hu and C. W. Shu, Weighted essentially non-oscillatory schemes on triangular meshes,
J. Comput. Phys.  150 (1999), 97-127.


\bibitem{ccp1} Peng Jin, Xi Deng and Feng Xiao, A direct ALE multi-moment finite volume scheme for the compressible Euler equations,
Commun. Comput. Phys.  24 (2018), 1300-1325.
\bibitem{ccp2} Yilang Liu, Weiwei Zhang and Chunna Li, A novel multi-dimensional limiter for high-order finite volume methods on unstructured grids,
Commun. Comput. Phys.  22 (2017), 1385-1412.

\bibitem{ccp3} Walter Boscheri and Raphael Loubere, High order accurate direct Arbitrary-Lagrangian-Eulerian ADER-MOOD finite volume schemes for non-conservative hyperbolic systems with stiff source terms, Commun. Comput. Phys.  21 (2017), 271-312.

\bibitem{ccp4}
Q. Wang, Y. X. Ren, J. H. Pan,  and W. A. Li, Compact high order finite volume method on unstructured grids III: Variational reconstruction, J. Comput. Phys. 337(2017), 1-26.

\bibitem{ccp5} Xiang Lai , Zhiqiang Sheng and Guangwei Yuan, Monotone finite volume scheme for three dimensional diffusion equation on tetrahedral meshes,
Commun. Comput. Phys.  21 (2017), 162-181.

\bibitem{toro04}V. A. Titarev and E. F. Toro, Finite-volume WENO schemes for
three-dimensional conservation laws,
J. Comput. Phys. 201 (2004), 238-260.

\bibitem{yuan2011}F. Teng, L. Yuan and T. Tang, A speed-up strategy for finite volume WENO schemes
for hyperbolic conservation laws,  Journal of Scientific Computing 46 (2011), 359-378.

\bibitem{improved16}P. Buchm{\"{u}}ller, J. Dreher and C. Helzel, Finite volume WENO methods for hyperbolic conservation laws on Cartesian grids with adaptive mesh refinement,
Applied Mathematics and Computation  272 (2016), 460-478.
\bibitem{improved18}P. Buchm{\"{u}}ller, J. Dreher and C. Helzel, Improved accuracy of high-order
WENO finite volume methods on Cartesian grids with adaptive mesh refinement,
Springer Proceedings in Mathematics \& Statistics 236 (2018), 263-272.

\bibitem{Tamaki2017CAF}Y. Tamaki and T. Imamura, Efficient dimension-by-dimension higher order finite-volume methods for a Cartesian grid with cell-based refinement, Comput. Fluids  144 (2017), 74-85.


\bibitem{wenoz}R. Borges, M. Carmona, B. Costa and W. S. Don, An improved weighted essentially non-oscillatory scheme for hyperbolic conservation laws,
J. Comput. Phys.  227 (2008), 3191-3211.


\bibitem{wenoz7}
W. S. Don, R. Borges,   Accuracy of the weighted essentially non-oscillatory conservative finite difference
schemes. J. Comput. Phys. 250 (2013), 347-372.


\bibitem{highweno}G. A. Gerolymos, D. S{\'{e}}n{\'{e}}chal and I. Vallet, Very-high-order {WENO} schemes,
J. Comput. Phys.  228 (2009), 8481-8524.
\bibitem{zhoutie01}T. Zhou, Y. F. Li and C. W. Shu, Numerical comparison of {WENO} finite volume and {Runge-Kutta} discontinuous {Galerkin} methods,
J. Sci. Comput. 16 (2001), 145-171.
\bibitem{bump1}J. O. Langseth and R. J. LeVeque, A wave propagation method for
three-dimensional hyperbolic conservation laws,
J. Comput. Phys.  165 (2000), 126-166.
\bibitem{tu2007}G. H. Tu and X. J. Yuan, A characteristic-based shock-capturing scheme for
hyperbolic problems,
J. Comput. Phys.  225 (2007), 2083-2097.
\bibitem{hushu02}J. Shi, C. Q. Hu and C. W. Shu, A technique of treating negative weights in {WENO} schemes,
J. Comput. Phys.  175 (2002), 108-127.
\bibitem{dengxiaogang17}Y. D. Dong, X. G. Deng, D. Xu and G. X. Wang, Reevaluation of high-order finite difference and finite volume algorithms with freestream preservation satisfied,
Comput. Fluids  156 (2002), 343-352.
\end{thebibliography}
\end{document}